\documentclass[12pt]{article}

\usepackage{amsmath}
\usepackage{amsthm}
\usepackage{eufrak}

\newtheorem{thm}{Theorem}

\newtheorem{lem}{Lemma}

\newtheorem*{cor*}{Corollary}

\newtheorem*{def*}{Definition}

\begin{document}

\title{A sharp inequality for the variance with respect to the Ewens Sampling Formula}

\author{\v Z. Baron\.{e}nas, E. Manstavi\v{c}ius and P. \v{S}apokait\.{e}}

\maketitle

\footnotetext{{\it AMS} 2000 {\it subject classification.} Primary
60C05;      secondary 05A16, 20P05. \break {\it Key words and
phrases}. Random permutation, linear statistics, quadratic form, eigenvalue,  discrete Hahn's polynomials, Tur\'an-Kubilius inequality.}

\begin{abstract}

The variance of a linear statistic defined  on the symmetric
  group endowed with the Ewens probability is
examined. Despite the dependence of the summands, it can be bounded
from above by  a constant multiple of the sum of variances of the summands. We find  the exact value of this
constant. The analysis of the appearing  quadratic forms and eigenvalue search is built upon the exponential matrices and discrete Hahn's polynomials.

\end{abstract}

\def\s{\smallskip}

\def\E{\mathbf{E}}
\def\C{\mathbf{C}}
\def\D{\mathbf{D}}
\def\G{\mathbf{G}}
\def\V{\mathbf{V}}
\def\Ra{\Rightarrow}
\def\N{\mathbf{N}}
\def\R{\mathbf{R}}
\def\S{\mathbf{S}_n}
\def\Z{\mathbf{Z}}
\def\k{\kappa}
\def\e{\varepsilon}
\def\ro{{\rm o}}

\def\rO{{\rm O}}
\def\re{{\rm e}}
\def\rd{{\rm d}}
\def\nt{\nu_{n,\theta}}
\def\Ent{{\mathbf E}_{n,\theta}}

\def\Env{{\mathbf E}_{n,1}}
\def\End{{\mathbf E}_{n,2}}

 \def\Varnt{{\mathbf{V}ar}_{n,\theta}}
 \def\Varnv{{\mathbf{V}ar}_{n,1}}
 \def\Varnd{{\mathbf{V}ar}_{n,2}}

\def\sj{\sum_{j\le n}}
\def\n{$n\to\infty$}
\def\cF{\mathcal F}

\section{Introduction and results}

The variance of a linear statistic defined  on the symmetric
  group endowed with the Ewens probability is
examined in the paper. The main obstacle to  overcome in this seemingly simple problem is the dependence of the summands. We propose an approach built upon  exponential matrices and special functions.

 Let $\S$ denote the symmetric group
 of permutations $\sigma$ acting  on $n\in\N$ letters. Each $\sigma\in\S$ has
 a unique    representation (up to the order) by the product of independent
cycles $\k_i$:
  \begin{equation}
           \sigma=\k_1\cdots\k_w                                                                                                                                                                                                                                                                                                                                                                         \label{1}
           \end{equation}
where $w=w(\sigma)$ denotes the number of cycles.
 Denote by $k_j(\sigma)\geq 0$  the number of cycles in (\ref{1}) of length $j$ for
 $1\leq j\leq n$  and introduce the \textit{cycle
vector} $\bar k(\sigma)=(k_1(\sigma),\dots,k_n(\sigma))$.

 As usual, set  $(x)_{m}=x(x+1)\cdots(x+m-1)$, $m\in\N_0:=\N\cup\{0\}$, for the increasing factorial. Denote also
\[
\Theta(m)=(\theta)_{m}/m!=[x^m](1-x)^{-\theta},
\]
 where  $[x^m] f(x)$ stands for the $m$th coefficient of a power series $f(x)$ and  $\theta>0$ is a parameter. The \textit{Ewens Probability Measure}  $\nt$  on $\S$ is defined by
\[
      \nt\big(\{\sigma\}\big)=\theta^{w(\sigma)}/(\theta)_{n}, \quad
      \sigma\in\S.
\]
Set $\ell(\bar s)=1s_1+\cdots +ns_n$ for a vector $\bar s=(s_1,\dots,s_n)\in {\N_0}^n$. Equality
 $\ell(\bar k(\sigma))=n$, valid for each $\sigma\in\S$, shows the dependence of the r.vs $k_j(\sigma)$ with respect to $\nu_{n,\theta}$. Throughout the paper, following the tradition of the probabilistic number theory, we leave out the elementary event $\sigma$ in all r.vs defined on the probability space $(\S, 2^{\S}, \nu_{n,\theta})$. It is well known (see, for example, \cite[Sect. 2.3]{ABT}) that the distribution of $\bar k(\sigma)$ can be written as the conditional distribution of $\bar\xi=(\xi_1,\dots,\xi_n)$, where $\xi_j$, $1\leq j\leq n$, are mutually independent Poisson
r.vs given on some probability space $\{\Omega, {\mathcal F}, {\mathbf P}\}$
with parameter ${\mathbf E}\xi_j=\theta/j$. Indeed,
\begin{equation}
   \nt\big(\bar k(\sigma)=\bar s\big)={\mathbf 1}\big\{\ell(\bar
   s)=n\big\}\Theta(n)^{-1}
     \prod_{j=1}^n\Big({\theta\over j}\Big)^{s_j}{1\over s_j!}={\mathbf P}\big(\bar\xi=\bar s\vert\,
     \ell(\bar\xi)=n\big).
     \label{ell}
\end{equation}
 Here ${\mathbf 1}\big\{\cdot\}$ stands for the indicator function. The
probability in (\ref{ell}), ascribed to the vector $\bar
s\in\N_0^n$, is called the {\itshape Ewens Sampling Formula}. It
has been introduced by W.~J.~Ewens \cite{Ew} to model the mutation
of genes.  For a comprehensive account of the recent applications of this ubiquitous distribution
in combinatorics and  statistics, see \cite{ABT}, \cite{J-K-B},  \cite{Feng}, or survey \cite{Crane} and the subsequent comments on it.

We prefer to stay within the theory of random permutations. Apart from $w(\sigma)$, other linear statistics (or {\textit{completely additive functions})
  \begin{equation}
  h(\sigma)\colon=a_1k_1(\sigma)+\cdots+ a_nk_n(\sigma),
 \label{h}
 \end{equation}
 where $\bar a:=(a_1,\dots, a_n)\in\R^n$ is a non-zero vector, continue to raise an interest. For example, $h(\sigma)$ with $a_j=\log j$, $j\leq
n$, is a good approximation for the logarithm of the
group-theoretical order of $\sigma\in\S$ (see \cite{ABT} or \cite{VZ-LMJ04}). The
case with  $a_j=\{x j\}$, where $\{u\}$ stands for the fractional
part of $u\in\R$, is met in the theory of random permutation
matrices (see \cite{Wieand}).

For an arbitrary $h(\sigma)$, the problem of finding  necessary and sufficient conditions, assuring the weak convergence of distributions
\begin{equation}
    \nt\big(h(\sigma)-\alpha(n)\leq x\beta(n)\big),
\label{Lthm}
\end{equation}
where $\alpha(n)\in\R$ and $\beta(n)\to\infty$ as $n\to\infty$, is still open (see \cite[Sect. 8.5]{ABT} or \cite{EM-OJM09} and the references therein). Obstacles  in the necessity part arise because of the dependence of the summands as shown by (\ref{ell}). This also happens  in the analysis of power moments carried out by the second
author \cite{EM-Prag06} and \cite{EM-AAM07} even in the
case $\theta=1$. Let us now focus on the variance.

By $\Ent f(\sigma)$ and $\Varnt f(\sigma)$ we denote the mean value
and the variance of a r.v. $f(\sigma)$ defined on $\S$ with
respect to $\nt$. For the particular function $h(\sigma)$ in
(\ref{h}), we also set   $A_{n,\theta}(\bar a)=\Ent h(\sigma)$ and
$
 D_{n,\theta}(\bar a) =\Varnt h(\sigma)$.
 Applying Watterson's  \cite{Wat} formulas (see \cite[(5.6), p. 96]{ABT}) for the factorial moments of $k_j(\sigma)$, one easily finds (see \cite{EMVS-AofA14} for the details) the expressions
\[
A_{n,\theta}(\bar a)=\theta\sum_{j\leq n} \frac{a_j}{j}\frac{\Theta(n-j)}{\Theta(n)}
\]
and
\begin{align}
D_n(\bar a)&=
\theta\sum_{j\leq n}{a_j^2\over j}\frac{\Theta(n-j)}{\Theta(n)}+\theta^2\sum_{i+j\leq n}\frac{a_ia_j}{i j}\frac{\Theta(n-i-j)}{\Theta(n)}\nonumber\\
&\quad -\theta^2\bigg(\sum_{j\leq n} \frac{a_j}{j}\frac{\Theta(n-j)}{\Theta(n)}\bigg)^2\nonumber\\
&=:
\theta B_n(\bar{a})+\theta^2\Delta_n(\bar{a}),
\label{DnB}
\end{align}
if $n\geq 2$ and
\[
B_n(\bar{a})=\sum_{j\leq n}{a_j^2\over j}\frac{\Theta(n-j)}{\Theta(n)}.
\]

The latter quantity is close to the sum of variances of the summands in the definition of $h(\sigma)$. In fact, formula (4) from
 \cite{EMVS-AofA14} shows that
 \[
    \sum_{j\leq n} a_j^2 \Varnt k_j(\sigma)-\theta B_n(\bar{a})=O\Big( n^{-(1\wedge\theta)}B_n(\bar{a})\Big),
    \]
    where $a\wedge b=\min\{a,b\}$ if $a,b\in\R$, with an absolute constant in  the symbol $O(\cdot)$.
We also have (see \cite{EMVS-AofA14})
\begin{equation}
 D_n(\bar{a})\leq C\theta B_n(\bar{a})
 \label{EMVS}
 \end{equation}
 uniformly in $n\geq 2$ with an absolute constant $C$ which can be specified.
  If $\theta\geq 1$, one can take  $C=2$.
   The purpose of the present paper is to find the exact value of $C$ in (\ref{EMVS}).

\begin{thm} \label{thm1} Let $\theta>0$ be arbitrary and $n\geq 2$. Then
\[
\tau_n(\theta):=\sup\bigg\{\frac{D_n(\bar a)}{\theta B_n(\bar{a})}: \quad \bar a\in\R^n\setminus\{0\}\bigg\}=\frac{\theta+2}{\theta+1}.
\]
The supremum is achieved taking $a_j=(\theta+2)j^2-(2n+\theta)j$ where $1\leq j\leq n$.
\end{thm}

The pioneering results  obtained in \cite{EM-Prag06} and \cite{EMZZ-LANA11} showed  that
$\tau_n(1)=3/2+O(n^{-1})$ and $\tau_n(2)=4/3+O(n^{-1})$. The approach originated in  Kubilius' paper \cite{JK} was based upon the extremal properties of the Jacobi polynomials. It unavoidably added a vanishing error term to the result.
 Recently  J.~Klimavi\v{c}ius and the second  author \cite{JKEM-AnBud18}  established that $\tau_n(1)= 3/2$ for all $n\geq 2$. Theorem 1 resumes the research for an arbitrary $\theta>0$. It is directly related to the above mentioned problem concerning  distributions (\ref{Lthm}). Applying  Theorem 1, we obtain that the weak convergence of (\ref{Lthm}) with $\beta(n)=\sqrt{\theta B_n(a)}$, which is natural to use, can take place  only to the limit laws having variance not exceeding $(\theta+2)/(\theta+1)$.

By virtue of (\ref{ell}), the result can be reformulated for the conditional variance of the linear statistics
\[
Y_n:=a_1\xi_1+\cdots+a_n\xi_n.
 \]
 We obtain  the following optimal inequality.

 \begin{cor*} Let $n\geq 2$ and $a_j\in{\mathbf R}$, $1\leq j\leq n$, be arbitrary. Then
 \[
     {\bf V}ar \big(Y_n\, \big| \ell(\bar\xi)=n\big)\leq\frac{\theta(\theta+2)}{\theta+1}
  \sum_{j=1}^n\frac{a_j^2}j\, \frac{\Theta(n-j)}{\Theta(n)}.
  \]
  \end{cor*}

The problem concerns the
quadratic forms $\Delta_n(\bar{a})$ and $B_n(\bar{a})$. The substitution
\[
a_j=\Big(\frac{j \Theta(n)}{\Theta(n-j)}\Big)^{1/2}\, x_j , \quad 1\leq j\leq
n,
\]
reduces  $B_n(\bar{a})$ to the square of Euclidean norm  $||\bar x||^2$ of the vector $\bar x=(x_1,\dots, x_n)\in\R^n$. Then $\Delta_n(\bar a)$ becomes a quadratic form, denoted afterwards by $\mathcal{M}_n(\bar x):=\bar x M_n {\bar x}'$, where ${\bar x}'$ is the column-vector and $M_n=((m_{ij}))$, $1\leq i,j\leq n$, is the matrix with  entries
 \begin{equation}
   m_{ij}=
    \frac{\Theta(n-i-j)}{\big(ij\Theta(n-i)\Theta(n-j)\big)^{1/2}}
    -\Big(\frac{\Theta(n-i)}{i\Theta(n)}\Big)^{1/2}\Big(\frac{\Theta(n-j)}{j\Theta(n)}\Big)^{1/2}.
      \label{mij}
      \end{equation}
Here we assume that $\Theta(-k)=0$ if $k\in\N$.  Now, by virtue of (\ref{DnB}),
\begin{align}
          \tau_{n,\theta}&=1+\theta\sup_{\bar
          x\not=\bar{0}}
          \Big(||\bar x||^{-2}{\mathcal M}_n(\bar x)\Big)=1+\theta\sup_{\bar x\not=\bar 0}\bigg(||\bar
          x||^{-2}\sum_{r=1}^{n}\mu_r x_r^2\bigg)\nonumber\\
          &=
          1+\theta\max_{1\leq r\leq n}\mu_r,
\label{min-eig}
\end{align}
 where $\{\mu_1,\dots, \mu_n\}$ is the spectrum of matrix $M_n$. So, Theorem 1 follows from the following proposition.

 \begin{thm} \label{ZB} The spectrum of the matrix $M_n$  comprises
 \[
 \mu_r=\frac{(-1)^{r}(r-1)!}{(\theta)_r}, \quad  1\leq r\leq n.
 \]
 For the eigenvector  corresponding to the maximal $\mu_2$, one may take the vector with coordinates
 \[
 \big((\theta+2)j-(2n+\theta)\big)(j\Theta(n-j))^{1/2}, \quad 1\leq j\leq n.
 \]
\end{thm}

The proof of Theorem \ref{ZB} presented in  the next section is built upon  exponential matrices.

The problem of  finding the remaining eigenvectors of  matrix $M_n$ also raises an interest. We solve it basing upon particular cases of the  generalized hypergeometric series which are exposed, for example, in  \cite[Chap. 5]{BW}. The hint to exploit them stems from  \cite{JKEM-AnBud18}. Let us confine ourselves to the case of  polynomials which, in the traditional notation, can be written as
\[
{}_{p+2}F_q\big(-m,-x,(a_p);  (b_q); z\big)=\sum_{k=0}^m \frac{(-m)_k (-x)_k (a_1)_k\cdots (a_p)_k}{(b_1)_k\cdots (b_q)_k k!} z^k,
\]
where $p,q, m\in \N_0$, $a_1,\dots, a_p; b_1,\dots, b_q\in \R$ are parameters.  Moreover, it suffices to reckon the discrete Hahn's polynomials
\[
    Q_r(x; \alpha,\beta, n)=\, {}_3F_2\big(-r, -x, r+\alpha+ \beta +1; \alpha+1, -n+1;1\big)
\]
by specifying  the parameters to $\alpha=1$ and $\beta=\theta-1$. In this case,
\[
q_r(x)= Q_r(x-1; 1,\theta-1, n),\quad  0\leq r\leq n-1,
\]
 have the following orthogonality property:
\begin{equation}
   <q_l,q_r>:=   \sum_{j=1}^n j q_l(j)q_r(j)\Theta(n-j)=\delta_{lr} \pi_r^2,
   \label{Inner}
   \end{equation}
   where $\delta_{lr}$ is the Kronecker symbol and $\pi_r>0$.
Note that,  up to a constant factor, $q_r(x)$, $1\leq r\leq n-1$, can be obtained uniquely  by the Gram--Schmidt orthogonalization procedure starting with $q_0(x)=1$ and form an orthogonal basis in the vector space of polynomials whose degrees do not exceed $n-1$. Exploiting this, in the isomorphic Euclidean space $\R^n$, we easily find the needed canonical basis for the matrix $M_n$.

 \begin{thm} \label{eigenv} The system of the vectors
\[
         \bar e_r=(e_{r1},\dots, e_{rn}), \quad 1\leq r\leq n,
         \]
 where
    \[
       e_{rj}=\pi_{r-1}^{-1}q_{r-1}(j)\sqrt{j\Theta(n-j)},\quad 1\leq j\leq n,
   \]
 is an orthonormal basis in $\R^n$. Moreover,  the vector $\bar e_r$ is the  eigenvector of  matrix $M_n$ corresponding to $\mu_r$ for each $1\leq r\leq n$.
 \end{thm}

The proof will be presented in the last section of the paper.

  Finally, the distributions of mappings defined on random permutations taken according to the Ewens probability are close to that defined on logarithmic decomposable combinatorial structures (see \cite{ABT}); therefore, we hope that our method is applicable when estimating the variances of similar statistics defined in such classes.

\section{Proof of Theorem 2}
 The idea is to find a matrix $L_n$ such that the product
\[
      {\rm e}^{L_n}M_n{\rm e}^{-L_n}=:\big(\big(w_{ij}\big)\big)
\]
is the triangle matrix with $w_{ij}=0$ if $1\leq j<i\leq n$ and
$w_{jj}=  \mu_j$ if $ 1\leq j\leq n$. This implies that the eigenvalues of $M_n$ are listed on the main diagonal of the product, as desired. The implementation has not been so evident.

At first,  we recall two identities from \cite{Prudnikov}.

\begin{lem} \label{lem1} Let $M, m\in\N_0$ and $ a,b\in\R$. Then
\begin{equation}
\sum_{k=0}^M{a+k\choose k}{b-k \choose M-k}=\sum_{k=0}^M{a+b-k\choose M-k}
\label{Pr1}
\end{equation}
and
\begin{equation}
\sum_{k=0}^M (-1)^k {M\choose k}{a-k \choose m}={a-M\choose m-M}.
\label{Pr2}
\end{equation}
\end{lem}

\textit{Proof}. See formulas (43) on page  618 and (56) on page 619 of \cite{Prudnikov}.

Let us introduce the matrix $L_n(\theta)=\big(\big(l_{ij}\big)\big)$ with the entries $l_{ij}=0$ for all $1\leq i,j\leq n$ but for $i=j+1$, where
\[
      l_{j+1,j}= -\bigg(\frac{(j+1)j \Theta(n-j-1)}{\Theta(n-j)}\bigg)^{1/2}, \quad 1\leq j\leq n-1.
      \]
Consider the powers
$
      L_n^k(\theta)=:\big(\big(l_{ij}^{(k)}\big)\big) $, $0\leq k\leq n-1$.
      The nonzero entries of $L_n^k(\theta)$ fill up the $k$th, $1\leq k\leq n-1$, diagonal under the main one. By induction, we observe that
\begin{align*}
                      l_{j+k,j}^{(k)}&=
                      l_{j+k,j+1}^{(k-1)} l_{j+1,j} \\
                      &=
                      l_{j+k,j+k-1} l_{j+k-1,j+k-2}\cdots l_{j+1,j}\\
                      &=
                      (-1)^k\prod_{r=0}^{k-1}\big( (j+r+1)(j+r)\big)^{1/2}\bigg(\prod_{r=0}^{k-1}\frac{\Theta(n-j-r-1)}{\Theta(n-j-r)}\bigg)^{1/2}\\
                      &=
                      (-1)^k (j)_k\Big(\frac{j+k}j\Big)^{1/2}\bigg(\frac{\Theta(n-j-k)}{\Theta(n-j)}\bigg)^{1/2}
                      \end{align*}
                      if $1\leq k\leq n-j$.
Hence the matrix
      $
       V:={\rm e}^{L_n(\theta)}=:\big(\big(v_{ij}\big)\big)
       $
has $v_{ij}=0$ if $1\leq i<j\leq n$ and
\begin{align*}
     v_{ij}&=\frac{l_{ij}^{(i-j)} }{(i-j)!}=         (-1)^{i-j}{i-1\choose j-1} \Big(\frac{i}{j}\Big)^{1/2}\bigg(\frac{\Theta(n-i)}{\Theta(n-j)}\bigg)^{1/2}\\
     &=
    (-1)^{i-j}{i\choose j} \Big(\frac{j}{i}\Big)^{1/2}\bigg(\frac{\Theta(n-i)}{\Theta(n-j)}\bigg)^{1/2}
      \end{align*}
    if $i\geq j$. Moreover, $V^{-1}={\rm e}^{-L_n(\theta)}=\big(\big(|v_{ij}|\big)\big)$ if $1\leq i,j\leq n$.

More technical obstacles arise  calculating
\begin{align*}
    w_{ij}&=\sum_{1\leq r\leq i\atop j\leq s\leq n} v_{ir}m_{rs}|v_{sj}| \\
    &=
    \sum_{1\leq r\leq i\atop j\leq s\leq n}(-1)^{i-r}{i\choose r} \Big(\frac{r}{i}\Big)^{1/2}\bigg(\frac{\Theta(n-i)}{\Theta(n-r)}\bigg)^{1/2}\cdot
    \frac{\Theta(n-r-s)}{\big(rs\Theta(n-r)\Theta(n-s)\big)^{1/2}}\\
    &\hskip 2 true in \times
    {s-1\choose j-1} \Big(\frac{s}{j}\Big)^{1/2}\bigg(\frac{\Theta(n-s)}{\Theta(n-j)}\bigg)^{1/2}\\
    &\quad -
    \sum_{1\leq r\leq i\atop j\leq s\leq n}(-1)^{i-r}{i\choose r} \Big(\frac{r}{i}\Big)^{1/2}\bigg(\frac{\Theta(n-i)}{\Theta(n-r)}\bigg)^{1/2}\cdot
    \Big(\frac{\Theta(n-r)}{r\Theta(n)}\Big)^{1/2}\Big(\frac{\Theta(n-s)}{s\Theta(n)}\Big)^{1/2}\\
    &\hskip 2true in \times
    {s-1\choose j-1} \Big(\frac{s}{j}\Big)^{1/2}\bigg(\frac{\Theta(n-s)}{\Theta(n-j)}\bigg)^{1/2}\\
    &=:\Sigma_1-\Sigma_2.
    \end{align*}
    Here
    \[
    \Sigma_1=(-1)^i\Big(\frac{\Theta(n-i)}{ij\Theta(n-j)}\Big)^{1/2}
    \sum_{1\leq r\leq i\wedge n-j}\frac{(-1)^{r}}{\Theta(n-r)}{i\choose r}\sum_{j\leq s\leq n-r}\Theta(n-r-s){s-1\choose j-1}.
    \]
After the change $s=n-r-k$, the inner sum reduces to that given in (\ref{Pr1}). So we obtain
\begin{align*}
&\sum_{k=0}^{n-r-j}{\theta-1+k\choose k}{n-r-1-k\choose n-r-j-k}
=
\sum_{k=0}^{n-r-j}{\theta+n-r-2-k\choose n-r-j-k}\\
&=
\sum_{l=0}^{n-r-j}{\theta-2+j+l\choose l}=[x^{n-r-j}]\frac1{(1-x)^{\theta+j}}
=
{\theta+n-r-1\choose n-r-j}.
\end{align*}
Hence
 \[
    \Sigma_1=(-1)^i\Big(\frac{\Theta(n-i)}{ij\Theta(n-j)}\Big)^{1/2}
    \sum_{1\leq r\leq i\wedge n-j}\frac{(-1)^{r}}{\Theta(n-r)}{i\choose r}{\theta+n-r-1\choose n-r-j}.
    \]
Similarly,
\[
    \Sigma_2=
    \frac{(-1)^i}{\Theta(n)}\Big(\frac{\Theta(n-i)}{ij\Theta(n-j)}\Big)^{1/2}\sum_{1\leq r\leq i}(-1)^{r}{i\choose r} \sum_{j\leq s\leq n}\Theta(n-s){s-1\choose j-1}.
    \]
Since
\[
\sum_{j\leq s\leq n}
    {s-1\choose j-1}\Theta(n-s)=[x^{n}]\bigg(\frac{x^j}{(1-x)^j}\cdot \frac1{(1-x)^\theta}\bigg)={\theta+n-1\choose n-j},
    \]
    we obtain
    \[
\Sigma_2=\frac{(-1)^{i+1}}{\Theta(n)}\Big(\frac{\Theta(n-i)}{ij\Theta(n-j)}\Big)^{1/2}{\theta+n-1\choose n-j}.
\]
Consequently,
\begin{align}
w_{ij}&=\Sigma_1-\Sigma_2\nonumber\\
&=
(-1)^i\Big(\frac{\Theta(n-i)}{ij\Theta(n-j)}\Big)^{1/2}
    \sum_{0\leq r\leq i\wedge n-j}\frac{(-1)^{r}}{\Theta(n-r)}{i\choose r}{\theta+n-r-1\choose n-r-j}\nonumber\\
    &=:
(-1)^i\Big(\frac{\Theta(n-i)}{ij\Theta(n-j)}\Big)^{1/2}\cdot \Sigma.
    \label{wij}
    \end{align}
Using the definition of $\Theta(m)$ and applying identity (\ref{Pr2}), we find that
\[
\Sigma=\frac{j!}{(\theta)_j}\sum_{0\leq r\leq i\wedge n-j}(-1)^{r}{i\choose r}{n-r\choose j}
=
\frac{j!}{(\theta)_j}{n-j\choose j-i}.
\]
Here $\Sigma=0$ if $j<i$ and $\Sigma=j!/(\theta)_j$ if $i=j$. Plugging this into (\ref{wij}), we obtain
\[
w_{ij}=(-1)^i\Big(\frac{\Theta(n-i)}{ij\Theta(n-j)}\Big)^{1/2}\frac{j!}{(\theta)_j}{n-j\choose j-i}=
\begin{cases} 0 &\quad \text{if}\; i>j,\\  (-1)^j(j-1)!/(\theta)_j &\quad \text{if}\; i=j.\end{cases}
\]

This proves the main assertion of
Theorem 2.

 It remains  to find
 \[
q_1(j)=
\big[(\theta+2)j-(2n+\theta)\big]/2(1-n),\quad 1\leq j\leq n,
 \]
  define $\bar e_2$, and  show that $\bar e_2M_n=\mu_2\bar e_2$. Since the latter is the subject of Theorem 3, we may omit the proof.

    Theorem 2 is proved.

    \smallskip

  {\bf Remark.}
The following observation is  worth mentioning. As it stems from the proof, $w_{ni}=w_{in}=\mu_n \delta_{in}$ if $1\leq i\leq n$, where $\delta_{in}$ is the Kronecker symbol. Hence the vector  $\bar v_n$ is the eigenvector of $M_n$ corresponding to $\mu_n$. Thus, it is proportional to the vector $\bar e_n$ from the next theorem.  To verify the mentioned property for $\bar e_n$, it suffice to apply a formula established even before the Gauss' seminal paper from 1812 on the hypergeometric series (see \cite[Lecture 7]{Askey} for the historical account). The Chu--Vandermonde formula ((7.16),  p.~59 of the same book) reads as follows:
\begin{equation}
{}_2F_1(-m, b; c;1)=\frac{(c-b)_m}{(c)_m}.
\label{Ch-V}
\end{equation}
  Hence
\[
     q_{n-1}(j)={}_2F_1(-j+1, \theta+n; 2;  1)= \frac{(2-n-\theta)_{j-1}}{(2)_{j-1}}=\frac{(-1)^{j-1} (\theta+n-j)_{j-1}}{j!} .
\]
This in turn yields that $\bar e_n M_n=\mu_n \bar e_n$.

\section{Proof of Theorem 3}

We now find all eigenvectors of the matrix $M_n$. Again, we have to recall a useful identity.

\begin{lem} \label{lem2} Let $p,q, M\in \N_0$, $\alpha,\beta\in\R$, and  $a_1, \dots, a_p; b_1\dots, b_q$ be the parameters such that the hypergeometric series below is correctly defined. Then
\begin{align*}
 &\sum_{k=0}^M {M\choose k} (\alpha)_{M-k} (\beta)_k  \cdot  {}_{p+1}F_q\big(-k, (a_p); (b_q); 1\big)\\
 &=
 (\alpha+\beta)_M \cdot {}_{p+2}F_{q+1}\big(-M, \beta,(a_p); \alpha+\beta, (b_q); 1\big).
\end{align*}
\end{lem}

\textit{Proof}. See formula (7) presented on page 388 in \cite{Prudn}.

As a corollary, we find the next sum.

\begin{lem} \label{lem3}  Let $n\geq 2$, $0\leq M\leq n-1$  and $0\leq r\leq n-1$. Then
\begin{align*}
   \Sigma_r(M)&:=\sum_{k=0}^M Q_r(k;1,\theta-1,n) \Theta(M-k)\\
   &=
   \frac{(\theta+1)_M}{M!} {}_4F_3\big(-M,1,-r, r+\theta+1; \theta+1, 2, 1-n; 1\big).
\end{align*}
\end{lem}

\textit{Proof}. Apply Lemma \ref{lem2}  for  $\alpha=\theta$, $\beta=1$, and $p=q=2$.

\smallskip

The obtained expressions of $\Sigma_{r-1}(M), 1\leq r\leq n$ will be used afterwards. For short, let
\[
{}_4F_3(-M)={}_4F_3(-M,1,1-r, r+\theta; \theta+1, 2, 1-n; 1).
\]

\begin{lem} \label{lem4} Let $\bar y_r=(y_{r1},\dots,y_{rn})=\pi_{r-1}\bar{e}_{r}M_n$, $1\leq r\leq n$, then
\begin{equation}
y_{ri}=-
\bigg(\frac{\Theta(n-i)}{i}\bigg)^{1/2} \frac{n}{r(r+\theta-1)}\cdot {}_3F_2\big(-r, -n+i, r+\theta-1; \theta, -n;1\big)
\label{yr}
\end{equation}
if $1\leq i\leq n$.
\end{lem}

\textit{Proof}. In the notation above,
\begin{align*}
  y_{ri}&=\frac1{(i\Theta(n-i))^{1/2}} \Sigma_{r-1}(n-i-1)-\bigg(\frac{\Theta(n-i)}{i}\bigg)^{1/2}\frac1{\Theta(n)} \Sigma_{r-1}(n-1)\\
  &=
\frac1{\theta}\bigg(\frac{\Theta(n-i)}{i}\bigg)^{1/2} \Big[(n-i)\cdot {}_4F_3(-n+i+1)
-n\cdot {}_4F_3(-n+1)\Big]
\end{align*}
if $1\leq i<n$. Since $(a)_{l-1}=(a-1)_l/(a-1)$ if $a\not=1$,  we have
\begin{align*}
 &  (n-i)\cdot{}_4F_3(-n+i+1)=(n-i)\sum_{l=1}^{r}\frac{(-n+i+1)_{l-1}(1-r)_{l-1}(r+\theta)_{l-1}}{(\theta+1)_{l-1}(-n+1)_{l-1} l!}\\
   &=
   -\frac{\theta n}{r(r+\theta-1)}\bigg[-1+\sum_{l=0}^{r}\frac{(-n+i)_l(-r)_l(r+\theta-1)_l}{(\theta)_l(-n)_l l!}\bigg]\\
   &=
    \frac{\theta n}{r(r+\theta-1)}\big[1-{}_3F_2\big(-r,-n+i, r+\theta-1; \theta, -n;1\big)\big].
   \end{align*}
   Similarly,
   \begin{align}
        n\cdot {}_4F_3(-n+1)&= \frac{\theta n}{r(r+\theta-1)}\big[1-{}_3F_2\big(-r, -n,  r+\theta-1; \theta, -n;1\big)\big]\nonumber\\
        &=
        \frac{\theta n}{r(r+\theta-1)}\big[1-{}_2F_1\big(-r, r+\theta-1; \theta;1\big)\big]\nonumber\\
        &=
               \frac{\theta n}{r(r+\theta-1)},
                \label{part}
                \end{align}
        by virtue of (\ref{Ch-V}). Hence
\begin{align*}
&(n-i)\cdot{}_4F_3(-n+i+1)-n\cdot{}_4F_3(-n+1)\\
&=
-\frac{\theta n}{r(r+\theta-1)}\cdot
{}_3F_2\big(-r, -n+i, r+\theta-1; \theta, -n;1\big).
\end{align*}
Plugging this into the previous expression of $y_{ri}$, we complete the proof in the case $i<n$.

If $i=n$, then, using Lemma 3 and (\ref{part}), we obtain
\begin{align*}
     y_{rn}&=-\frac1{\sqrt{n} \Theta(n)} \Sigma_{r-1}(n-1)= -\frac1{\sqrt{n} \Theta(n)}\cdot \frac{(\theta+1)_{n-1}}{(n-1)!} \cdot \frac{\theta}{r(r+\theta-1)}\\
     &=
      -\frac{\sqrt{n}}{r(r+\theta-1)}.
      \end{align*}
      This is consistent with expression (\ref{yr}) given in the  lemma.

  Lemma \ref{lem4} is proved.

\medskip
\textit{Proof of Theorem} 3. Let $1\leq r\leq n$ be fixed. Recall that  $q_k(x), 0\leq k\leq r$, span the subspace of polynomials whose degrees do not exceed $r$. Analyse the polynomial appearing in Lemma \ref{lem4}, namely,
\[
    \Phi_r(x)=\, {}_3F_2\big(-r, -n+x, r+\theta-1; \theta, -n;1\big).
    \]
 As we have seen proving  (\ref{part}), we have $ \Phi_r(0)=0$. Hence $ x^{-1}\Phi_r(x)$ is a polynomial of degree $r-1$ and
\[
 \Phi_r(x)=x\sum_{k=0}^{r-1} c_k q_r(x), \quad c_k\in \R.
 \]
  The leading coefficients of the polynomials $\Phi_r(x)$ and $q_r(x)$ are, respectively,
\[
\frac{(-1)^r(r+\theta-1)_r}{(\theta)_r (-n)_r }, \qquad \frac{(r+\theta)_{r-1}}{r! (-n+1)_{r-1}}.
\]
Consequently,
\[
   c_{r-1}=\frac{(-1)^r(r+\theta-1)_r}{(\theta)_r (-n)_r }\cdot \frac{r!(-n+1)_{r-1}}{(r+\theta)_{r-1}}=
   \frac{(-1)^{r-1} r!(r+\theta-1)}{(\theta)_r n }.
   \]
   Now, the result of Lemma \ref{lem4} can be rewritten as follows:
   \begin{align*}
   y_{ri}&=-\bigg(\frac{\Theta(n-i)}{i}\bigg)^{1/2} \frac{n i}{r(r+\theta-1)}\bigg[c_{r-1}q_{r-1}(i)+\sum_{k=0}^{r-2} c_k q_k(i)\bigg]\\
   &=
   \big(i \Theta(n-i)\big)^{1/2}\bigg[\frac{(-1)^r (r-1)!}{(\theta)_{r}} q_{r-1}(i)+\sum_{k=0}^{r-2} d_k q_k(i)\bigg]
   \end{align*}
   with some coefficients $d_k=d_k(n,r,\theta)$ for each $1\leq i,r\leq n$. Note that the fraction in the brackets is just $\mu_r$ found in Theorem 2.

For  $1\leq r\leq l\leq n$, applying the last formula and the definition of the inner product (\ref{Inner}), we obtain
\begin{align*}
   \bar e_lM_n \bar e_r'&= \pi_{r-1}^{-1}\bar e_l y_r' = \frac{1}{\pi_{l-1}\pi_{r-1}}\bigg[\mu_r<q_{l-1},q_{r-1}>+
   \sum_{k=0}^{r-2}d_k <q_{l-1},q_k>\bigg]\\
   &=\mu_r\delta_{rl}
   \end{align*}
   by virtue of orthogonality.
  This shows that each $\bar e_r$ in the basis is the eigenvector for $M_n$ corresponding to $\mu_r$.

  Theorem 3 is proved.

\medskip

{\bf Concluding remark.} Comparing the expression $\bar e_r M_n=\mu_r\bar e_r$ with $\pi_{r-1}^{-1} \bar y_r$ given by Lemma \ref{lem4}, we arrive to a seemingly new relation of the generalized hypergeometric functions. For $1\leq i,r\leq n$ and $\theta>0$, it holds that
\begin{align*}
   (-1)^{r-1} r! i \cdot {}_3F_2\big(-r+1,& -i+1, r+\theta; 2, -n+1;1\big)\\
   &=
   (\theta)_{r-1} n \cdot {}_3F_2\big(-r, -n+i, r+\theta-1; \theta, -n;1\big).
   \end{align*}
Derivation of it using an appropriate  sequence of the so-called contiguous relations (see \cite{Rainv})
 would not be short.

\vskip 0.5 true in
Affiliation of the authors : Institute of Mathematics, Vilnius University, Naugarduko str. 24, LT-03225 Vilnius, Lithuania

Corresponding author: Eugenijus Manstavi\v{c}ius, 

 email: eugenijus.manstavicius@mif.vu.lt


\begin{thebibliography}{9}

 \bibitem{ABT} Arratia,~R., Barbour,~A.\,D., and Tavar\'e,~S. (2003) {\em Logarithmic
        Combinatorial Structures: a Probabilistic Approach}. EMS
        Monographs in Mathematics, EMS Publishing House, Z\"urich.

 \bibitem{Askey} Askey,~R. (1975) {\em Polynomials and Special Functions}. Soc. Industr. and Applied Maths, Philadelphia, Pensylvania.

\bibitem{BW} Beals,~R. and Wong,~R. (2012) {\em Special Functions}. Cambridge University Press, Cambridge.

\bibitem{Crane} Crane,~H. (2016) The ubiquitous Ewens sampling formula. {\em Statist. Sci.} {\bf 31} 1--19.

 \bibitem{Ew} Ewens,~W.\,J. (1972)  The sampling theory of selectively
            neutral alleles. {\em Theor.  Pop.  Biol.}  {\bf 3} 87--112.


\bibitem{Feng} Feng,~Sh. (2010) {\em The Poisson-Dirichlet Distribution and Related Topics}. Springer, Berlin.

  \bibitem{J-K-B} Johnson,~N.\,S., Kotz,~S., and Balakrishnan,~N. (1997) {\em Discrete Multivariate Distributions}.
   Wiley, New York.

\bibitem{JKEM-AnBud18}  Klimavi\v{c}ius,~J., and Manstavi\v cius,~E. (2018) The Tur\'an-Kubilius inequality on permutations. {\em Annales Univ. Sci. Budapest., Sect. Comp.} {\bf48} 45--51.

\bibitem{JK} Kubilius,~J. (1983) Estimating the second central
moment for strongly additive arithmetic functions (Russian). {\em Liet.
matem. rink.},  {\bf 23} 110--117; translation in {\em Lith. Math. J.} {\bf23}  61--69.

  \bibitem{EM-Prag06} Manstavi\v cius,~E. (2006) Conditional Probabilities in Combinatorics. The Cost of Dependence.
In: {\em Prague Stochastics}. 7th Prague symposium on
asymptotic statistics and 15th Prague conference on information
theory, statistical decision functions and random processes (M.~Hu\v skov'a and M. Jan\v zura, eds), Prague, Charles University, Matfyzpress,  523--532.


\bibitem{EM-AAM07} Manstavi\v cius,~E. (2007) Moments of additive functions defined on the symmetric group,
{\em Acta Appl. Math.} {\bf 97} 119--127.

\bibitem{EM-OJM09} Manstavi\v cius,~E. (2009) An analytic method in probabilistic combinatorics. {\em Osaka J. Math.} {\bf46}  273–--290.

\bibitem{EMVS-AofA14} Manstavi\v{c}ius,~E., and Stepanauskas,~V. (2014)  On variance of an additive function with respect to a generalized Ewens probability. In: {\em Proceedings of the 25th International Conference on Probabilistic, Combinatorial and Asymptotic Methods for the Analysis of Algorithms, Discrete Math. Theor. Comput. Sci. Proc.}, BA, Assoc. Discrete Math. Theor. Comput. Sci., Nancy, 301–--311.

\bibitem{EMZZ-LANA11} Manstavi\v{c}ius,~E., and \v{Z}ilinskas,~\v{Z}. (2011)  On a variance related to the Ewens sampling formula. {\em Nonlinear Anal. Model. Control} {\bf16} 453--466.

\bibitem{Prudnikov} Prudnikov,~A.\,P., Brychkov,~Yu.\,A., and Marichev,~O.\,I. (1998) {\em Integrals and Series}, vol. 1. Taylor \& Francis, Fourth printing.

\bibitem{Prudn} Prudnikov,~A.\,P., Brychkov,~Yu.\,A., and Marichev,~O.\,I. (1990) {\em Integrals and Series}, vol. 3, More Special Functions. Gordon and Breach Sci. Publ., New York.


\bibitem{Rainv} Rainville,~E.\,D. (1960) {\em Special Functions}. Macmillan, New York.


  \bibitem{Wat} Watterson,~G.\,A. (1974) The sampling theory of selectively
neutral alleles diffusion model. {\em Advances in Appl. Probab.},
{\bf6} 463--488.

\bibitem{Wieand} Wieand,~K. (1998) {\em Eigenvalue Distributions of
         Random Matrices in the Permutation Group and Compact Lie
         Groups}. PhD thesis, Harvard University.

\bibitem{VZ-LMJ04} Zacharovas,~V. (2004) Distribution of the logarithm of the order of a random permutation (Russian). {\em Liet. mat. rink.} {\bf 44} 372--406; translation in {\em Lithuanian Math. J.} {\bf44} 296–--327.
\end{thebibliography}
\end{document}